\documentstyle[12pt]{article}
\newcommand{\ba}{\begin{array}}
\newcommand{\ea}{\end{array}}
\newcommand{\be}{\begin{equation}}
\newcommand{\ee}{\end{equation}}
\newcommand{\ben}{\begin{enumerate}}
\newcommand{\een}{\end{enumerate}}

\newcommand{\flip}{\tau}
\newcommand{\Poly}{\mbox{\rm Poly}\, }

\newtheorem{Thm}{Theorem}[section]
\newtheorem{Def}[Thm]{Definition}
\newtheorem{Prop}[Thm]{Proposition}

\newtheorem{Rem}[Thm]{Remark}

\newcommand{\ov}{\overline}
\newcommand{\cj}{{\ov{1}}}
\newcommand{\cd}{{\ov{2}}}
\newcommand{\cu}{\ov{u}}
\newcommand{\cA}{\ov{A}}
\newcommand{\cB}{\ov{B}}

\newcommand{\cD}{\ov{D}}
\newcommand{\cL}{\ov{L}}

\newcommand{\Mor}{\mbox{\rm Mor}\, }
\newcommand{\id}{\mbox{\rm id}\, }

\font \msb=msbm10 scaled \magstep1

\newcommand{\bR}{\mbox{\msb R} }
\newcommand{\bC}{\mbox{\msb C} }

\newcommand{\ar}{\alpha }
\newcommand{\br}{\beta }
\newcommand{\gr}{\gamma }
\newcommand{\dr}{\delta }
\newcommand{\er}{\varepsilon }

\def\Dr{\Delta }
\newcommand{\cP}{{\cal P}}
\newcommand{\cQ}{{\cal Q}}
\newcommand{\pA}{{\cal A}}
\newcommand{\tpA}{{\widetilde{\cal A}}}
\newcommand{\tis}{{\widetilde{s}}}
\newcommand{\pC}{{\cal C}}
\newcommand{\tpC}{{\widetilde{\cal C}}}
\newcommand{\pJ}{{\cal J}}
\newcommand{\lel}{\left\langle}
\newcommand{\rr}{\right\rangle}

\begin{document}

\title{\bf Quantum Lorentz and braided Poincar\'{e}
groups\thanks{Supported by Polish KBN grant No 2 P301 020 07}}
\author{{\bf S. Zakrzewski}  \\
\small{Department of Mathematical Methods in Physics,
University of Warsaw} \\ \small{Ho\.{z}a 74, 00-682 Warsaw, Poland}}
\date{}
\maketitle

\begin{abstract}
Quantum Lorentz groups $H$ admitting quantum Minkowski space $V$
are selected. Natural structure of a quantum space $G=V\times H$
is introduced, defining a quantum group structure on $G$ only for
triangular $H$ ($q=1$). We show that it defines a braided quantum group
structure on $G$ for $|q|=1$.
\end{abstract}

\section{Introduction}

Any example of a quantum Poincar\'{e} group \cite{PP} is
constructed using one of quantum Lorentz groups introduced in
\cite{WZ}. However, only very special cases of the latter
(triangular deformations) can be used for this purpose. Cases
related to the celebrated $q$-deformation of Drinfeld and Jimbo
are, unfortunately, excluded.  This is in fact a general feature
of inhomogeneous quantum groups \cite{PP1,pspg}.

It turned out recently that this obstacle can be circumvented,
if one allows the deformed inhomogeneous group to be a braided
quantum group rather than an ordinary quantum group. It means
that the comultiplication is a morphism into a nontrivial
crossed-product algebra rather than the usual product. It turns
out that on the level of generators, the only non-trivial
cross-relations are those for the translation coordinates.
These results have been derived in \cite{cb} for the case when the
homogeneous part is the standard $q$-deformed (with $|q|=1$)
orthogonal quantum group $SO(p,p)$, $SO(p,p+1)$ \cite{FRT} or
$SO(p,p+2)$ \cite{it}. I have learned recently about the paper
\cite{Drab} where results of similar type (without the
reality condition) were obtained (cf. also \cite{Ma,Ma1}).

In the present paper we study the case when the homogeneous part
$H$ is the Lorentz group. This case requires a separate study,
because we have the possibility to take into account the
complete classification of quantum deformations \cite{WZ}.
Another reason for a separate treatment is that we want to
consider the `more fundamental' simply connected
$SL(2,\bC )$ group instead of $SO(1,3)$.

The paper is organized as follows. In Section 1 we recall
non-triangular, deformation-type cases of quantum Lorentz group
$H$. In Section 2 we select those cases which have the corresponding
quantum Minkowski space $V$ (this happens for $|q|=1$ or $q^2\in
\bR$). In Section 3 we construct a natural crossed `cartesian product' $G$
of $V$ and $H$ (as quantum spaces). In Section 4 we
investigate conditions under which the natural formula for the
comultiplication on generators, defines a morphism of algebras,
the product algebra being understood with suitable crossed (or
braided) structure.

The same program on the Poisson level has already been presented in
\cite{cb}.

We conclude in Section 5 with explicit commutation relations for
the Minkowski space. Several proofs are shifted to
the Appendix.

\section{Quantum Lorentz groups}

We recall that the $^*$-algebra ${\cal A} = \Poly (H)$ of
polynomials on quantum $H=SL(2,\bC )$ is generated by matrix
elements of
$$u = (u^A_B)_{A,B=1,2} =\left( \ba{cc} u^1_1 & u^1_2\\ u^2_1 &
u^2_2\ea\right)$$
subject to relations
\be\label{uuE}
u_1u_2E=E,\qquad E'u_1u_2=E',\qquad Xu_1\ov{u}_2=\ov{u}_1u_2X,
\ee
where $E$, $E'$ and $X$ are described in Theorem~2.2 of \cite{WZ}.
Here the subscripts 1 and 2 refer to the position of a given
object in the tensor product of the underlying `arithmetic'
vector space (in this case $\bC ^2$, with the standard basis
$e_1, e_2$). For instance, the first equality means that
$u^A_Cu^B_DE^{CD}=E^{AB}$ (summation convention). We omit the
subscripts when the object has only one natural position in a
given situation (like $E$ for instance). The complex conjugate
$\cu$ of $u$ is given by
$$ \cu = (\cu ^{\cA}_{\cB})_{A,B=1,2} =\left( \ba{cc}(u^1_1)^* &
(u^1_2)^*\\ (u^2_1)^* & (u^2_2)^*\ea\right),\qquad \mbox{i.e.}
\;\;\;\cu ^{\cA}_{\cB} =(u^A_B)^* .$$
The `barred' indices refer to the complex conjugated basis
$e_\cj :=\ov{e_1}$, $e_\cd :=\ov{e_2}$ in $\ov{\bC ^2}$.

With the standard comultiplication defined on generators by
$\Dr u=uu'$ (primed coordinates refer to the {\em second copy}
of $H$; in a less compact notation, $\Dr u^A_B=u^A_C\otimes
u^C_B\in {\cal A}\otimes {\cal A}$), the above $^*$-algebra
becomes a Hopf $^*$-algebra. ($\Dr$ preserves the relations, for
instance $\Dr u_1 \Dr u_2 E = u_1u'_1u_2u'_2E =
u_1u_2u'_1u'_2E = u_1u_2E=E$.)

In the sequel we focus on non-triangular deformations.
It means that
\be
  E=e_1\otimes e_2-qe_2\otimes e_1,\qquad E'=e^2\otimes
e^1-q^{-1}e^1\otimes e^2\qquad q\in \bC \setminus \{ 0,i,-i\}
\ee
(the standard $q$-deformation) and $X$ is given by (13) or (15)
of \cite{WZ}, i.e. we have one of the following two cases
\ben
\item
$X=t^{\frac12}(e^1_{\ov{1}}\otimes e^{\ov{1}}_1+
e^2_{\ov{2}}\otimes e^{\ov{2}}_2)+t^{-\frac12}
(e^1_{\ov{2}}\otimes e^{\ov{2}}_1+e^2_{\ov{1}}\otimes e^{\ov{1}}_2)$
\ \ \ \ \ \ $0<t\in \bR$
\item
$X=q^{\frac12}(e^1_{\ov{1}}\otimes e^{\ov{1}}_1+
e^2_{\ov{2}}\otimes e^{\ov{2}}_2)+
q^{-\frac12}(e^1_{\ov{2}}\otimes e^{\ov{2}}_1+e^2_{\ov{1}}\otimes
e^{\ov{1}}_2)\pm q^{\frac12}e^2_{\ov{1}}\otimes e^{\ov{2}}_1 $ \hspace{1cm}
(for $0<q\in \bR$)
\een
(with the obvious notation for the matrix units
$e^1_{\cj}:=e_{\cj}\otimes e^1$ etc.).

Any matrix which intertwines $u_1u_2$ with itself and satisfies
the braid equation is proportional to
\be\label{M}
M:= q\cP '-q^{-1}\cP\qquad \mbox{or}\qquad M^{-1}=
q^{-1}\cP '-q\cP ,
\ee
where $\cP := - (q+q^{-1})^{-1}EE'$ (the deformed
antisymmetrizer) is the projection on $E$ parallel to $\ker E'$
and $\cP ':= I-\cP $ (the deformed symmetrizer).
Conjugating $M^{\pm 1}u_1u_2=u_1u_2M^{\pm 1}$
we obtain $K^{\pm 1}\ov{u}_1\ov{u}_2=\ov{u}_1\ov{u}_2K^{\pm 1}$,
where
\be\label{K}
K := \flip \ov{M}\flip  = \ov{q}\cQ '-\ov{q}^{-1}\cQ ,\qquad
\cQ :=\flip \ov{\cP}\flip,\;\;\; \cQ ' :=I-\cQ  .
\ee
Throughout the paper, $\flip$ denotes the permutation in the
tensor product.

\section{Quantum Minkowski spaces}

In order to discuss quantum Minkowski spaces, covariant under
the quantum Lorentz group, we consider the four-dimensional
representation of the latter,
\be
h:= u_1\cu_2 \qquad \mbox{i.e.}\;\;\;
h^{A\cB }_{C\cD }:=u^A_C\cu ^{\cB}_{\cD}.
\ee
Note that $\flip \ov{h}\flip =h$. It means that in a basis
of elements selfadjoint with respect to the natural conjugation
$x\mapsto \flip \ov{x}$ in $\bC ^2\otimes \ov{\bC ^2}$, such as the
basis of Pauli matrices $\sigma ^{A\cB}_j$ ($j=0,\ldots,3$), the
matrix $h$ has selfadjoint elements. In considerations which
refer only to the four-dimensional (`vector') representation, it
is often convenient to use exactly the components of $h$ in the
basis of Pauli matrices. These components will be denoted by
$h^j_k$ ($j,k=0,\ldots ,3$). In the leg-numbering notation (like
in (\ref{uuE})) we shall use bold subscripts for the
four-dimensional case. For instance, the tensor square of $h$
will be denoted either by $h_{12}h_{34}$ (referring to the spinor
representation) or by $h_{\bf 1}h_{\bf 2}$ (referring to the
vector representation).

Now we look for appropriate quadratic commutation relations
defining the quantum Minkowski space. We use here the standard
method of dealing with `quantum vector spaces'. The algebra of
polynomials on quantum Minkowski space should be generated by
four generators $x = (x^{A\cB })_{A,B=1,2} = (x^j)_{j=0,\ldots ,3}$
satisfying the reality condition
\be
\flip \ov{x} = x\qquad (\mbox{i.e.}\;\; (x^{A\cB})^*=x^{B\cA}\;\;\;
\mbox{or}\;\;\; (x^j)^*=x^j)
\ee
and some quadratic relations $Ax_{\bf 1}x_{\bf 2}=0$, such that
$\Dr_V x:=hx'$ satisfies the same relations (note, that $\Dr_Vx$
satisfies the reality automatically). The last requirement will
be satisfied if $A$ is an intertwiner of $h_{\bf 1}h_{\bf 2}$:
\be
A\Dr_V x_{\bf 1}\Dr_V x_{\bf 2}=Ah_{\bf 1}x_{\bf 1}h_{\bf 2}x_{\bf 2}=
Ah_{\bf 1}h_{\bf 2}x_{\bf 1}x_{\bf 2}=h_{\bf 1}h_{\bf 2}Ax_{\bf
1}x_{\bf 2}=0
\ee
(this is the key point of the method of \cite{FRT}). It remains
to choose an appropriate intertwiner: it should be a deformation
of the antisymmetrizer (see also Remark~\ref{rema} below).

{}From $M^{\pm 1}$, $K^{\pm 1}$, $X$ and $X^{-1}$ we can build
easily four intertwiners of $h_{\bf 1}h_{\bf 2}=h_{12}h_{34}$, namely
\be\label{Wpm}
\hat{R}_{\pm} : = X_{23}(M_{12}K^{\pm 1}_{34})X^{-1}_{23}
\ee
and their inverses. Each of them becomes the permutation in
the classical limit.
\begin{Prop}\label{braid}
Matrices $\hat{R}_{\pm}$ satisfy the braid equation (with $\bC
^2\otimes \ov{\bC ^2}$ being the elementary space).
\end{Prop}
For the proof we refer to Appendix (Section~\ref{biM}).

Substituting (\ref{M}), (\ref{K}) into (\ref{Wpm}) we get the
spectral decomposition
\be
\hat{R}_{\pm} = X_{23}(q\ov{q}^{\pm 1}\cP '\otimes \cQ ' +
q^{-1}\ov{q}^{\mp 1 }\cP\otimes \cQ -
q\ov{q}^{\mp 1} \cP '\otimes \cQ - q^{-1}\ov{q}^{\pm 1}\cP\otimes
\cQ ')X^{-1}_{23}.
\ee
Since the projections $\cP '\otimes \cQ$, $\cP\otimes\cQ '$ are
3-dimensional,
\be
P^{(-)}:= X_{23}(\cP '\otimes \cQ + \cP\otimes \cQ ' )X^{-1}_{23}.
\ee
is a good candidate for the deformed antisymmetrizer. It is in
fact easy to see that it becomes the classical antisymmetrizer
in the classical limit.

\begin{Rem}\label{rema}
{\rm
It is not necessary to use the argument of a `deformed
antisymmetrizer'. In fact, there is a more straightforward
(logical) approach. Note that the subspace $V^*$ spanned by $x^j$ is
invariant with respect to $H$, and we are looking just for a
6-dimensional invariant subspace of $V^*\otimes V^*$. It must be
therefore the direct sum of the two 3-dimensional
irreducible subrepresentations in $V^*\otimes V^*$.
}
\end{Rem}

\begin{Def}
{\rm The $^*$-algebra generated by $(x^j)^*=x^j$ (i.e.
$(x^{A\ov{B}})^* = x^{B\ov{A}}$) and relations
\be\label{Mink}
P^{(-)}x_{\bf 1}x_{\bf 2} =0 \qquad (\mbox{i.e.} \;\; P^{(-)}x_{12}x_{34}=0)
\ee
is said to be the} $^*$-algebra of polynomials on quantum
Minkowski space {\rm (and denoted by $\Poly (V)$) if it has the
classical size (i.e. if the Poincar\'{e}-Birkhoff-Witt theorem
holds)}.
\end{Def}
\begin{Prop}
Quantum Minkowski spaces exist only for
\be\label{warun}
|q|=1\qquad \mbox{or} \qquad \ov{q}^2=q^2.
\ee
\end{Prop}
For the proof, see Appendix (Section~\ref{spM}).
Note that this result was conjectured in \cite{stand}.

In the sequel we assume one of the two possibilities: $q=\ov{q}$
or $|q|=1$ (we discard $q=-\ov{q}$ as not of deformation type).

Note that
\begin{eqnarray}
\mbox{for}\;\; q=\ov{q},\qquad
\hat{R}_{+} & = & X_{23}(q^2\cP '\otimes \cQ '
+ q^{-2}\cP\otimes \cQ - \cP '\otimes \cQ -
 \cP\otimes \cQ ')X^{-1}_{23} \\
\hat{R}_{-} & = & X_{23}(\cP '\otimes \cQ '
+ \cP\otimes \cQ - q^2\cP '\otimes \cQ -
q^{-2}\cP\otimes \cQ ')X^{-1}_{23} \label{R-1}\\
\mbox{for}\;\; |q|=1,\qquad
\hat{R}_{+} & = & X_{23}(\cP '\otimes \cQ '
+ \cP\otimes \cQ - q^2\cP '\otimes \cQ -
q^{-2}\cP\otimes \cQ ')X^{-1}_{23} \\
\hat{R}_{-} & = & X_{23}(q^2\cP '\otimes \cQ '
+ q^{-2}\cP\otimes \cQ - \cP '\otimes \cQ -
 \cP\otimes\cQ ' )X^{-1}_{23}, \label{R-2}
\end{eqnarray}
hence
\begin{eqnarray*}
\mbox{for}\;\; q=\ov{q},\qquad
\mbox{(\ref{Mink})} & \Longleftrightarrow &
\hat{R}_- x_{\bf 1}x_{\bf 2}=x_{\bf 1}x_{\bf 2} , \\
\mbox{for}\;\; |q|=1,\qquad
\mbox{(\ref{Mink})} & \Longleftrightarrow &
\hat{R}_+ x_{\bf 1}x_{\bf 2}=x_{\bf 1}x_{\bf 2},
\end{eqnarray*}
which `explains' why for $q=\ov{q}$ or $|q|=1$ we get the
appropriate size of the algebra generated by $x$, namely,
different ways of ordering the polynomials of the third degree
give the same result, due to the Yang-Baxter property of
$\hat{R}_{\pm}$ (Prop.~\ref{braid}):
\be\label{YBa}
R_{\bf 12}R_{\bf 13}R_{\bf 23}x_{\bf 1}x_{\bf 2}x_{\bf 3}  =
x_{\bf 3}x_{\bf 2}x_{\bf 1} =
R_{\bf 23}R_{\bf 13}R_{\bf 12}x_{\bf 1}x_{\bf 2}x_{\bf 3},
\ee
where $R=\tau \hat{R}_-$ for $\ov{q}=q$ and $R=\tau\hat{R}_+$
for $|q|=1$.

\section{Crossed product of Minkowski with Lorentz}

In this section we shall introduce a crossed tensor product of
$\Poly (H)$ and $\Poly (V)$ in such a way, that the standard
comultiplication
\be
\Dr u = uu',\qquad \Dr x = x+hx',
\ee
preserves `as much as possible' from the algebraic structure
(preserves as many relations as possible). Technically (see
Theorem~\ref{Th} below for the precise statement), we consider
the universal $^*$-algebra ${\cal B}$ generated by $u^A_B$ and
$x^j=(x^j)^*$, satisfying (\ref{uuE}), (\ref{Mink}) and the {\em
cross relations}
\be\label{cross}
x_{12}u_3 = Tu_1x_{23}\qquad (\mbox{i.e.}\;\; x^{A\cB}u^C_D=
T^{A\cB C}_{EK\cL}u^E_Dx^{K\cL})
\ee
for an appropriate matrix $T$, which we select after some
discussion.

Note that the `preservation of relations' by $\Dr$ means that
$\Dr u$ and $\Dr x$ do satisfy the same relations as $u$ and
$x$. Let us check when it happens. Of course, $\Dr u$ satisfies
(\ref{uuE}) as before. Since
$$
\Dr x_{12}\Dr u_3 = (x_{12}+h_{12}x'_{12})u_3u'_3 =
 Tu_1x_{23}u'_1 + h_{12}u_3Tu'_1x'_{23}
$$
(we use here $x'_{12}u_3=u_3x'_{12}$) and
$$
T\Dr u_1 \Dr x_{12}=T u_1u'_1(x_{23}+h_{23}x'_{23})=
 Tu_1x_{23}u'_1 + Tu_1h_{23}u'_1x'_{23},
$$
$\Dr x$ and $\Dr u$ satisfy (\ref{cross}) if
\be
Tu_1h_{23}=h_{12}u_3T,
\ee
i.e. $T\in \Mor (u_1u_2\cu _3,u_1\cu _2 u_3)$ ($T$ intertwines
$u_1u_2\cu_3$ with $u_1\cu_2u_3$). It means that
\be
T=X_{23}S_{12}
\ee
for some $S\in \Mor (u_1u_2,u_1u_2)$, which we assume to be
invertible.

The discussion of when $\Dr$ preserves (\ref{Mink}) will be
postponed till the next section.

By taking the star operation of (\ref{cross}), we obtain
\be
\cu_1x_{23}=X_{12}(\tau \ov{S}\tau )_{23}x_{12}\cu_3
\ee
(we have used the property $\tau\ov{X}\tau =X$), hence
\be\label{xcu}
x_{12}\cu_3 = (\tau \ov{S}^{-1}\tau)_{23}X^{-1}_{12}\cu_1x_{23}.
\ee
It follows that
\be\label{xh}
x_{12}h_{34}=X_{23}S_{12}(\tau
\ov{S}^{-1}\tau)_{34}X^{-1}_{23}h_{12}x_{34}
\ee
and the matrix governing the commutation of $x$ and $h$ has
similar structure to $\hat{R}_-$ in (\ref{Wpm}). It suggests that $S$
should be proportional to $M$ or $M^{-1}$. We shall show that it
is indeed the case, if we require ${\cal B}$ to have a correct size.

Recall that a {\em crossed tensor product} of two
algebras, ${\cal C}$ and ${\cal D}$, is the tensor product of
vector spaces ${\cal C}\otimes {\cal D}$ equipped with the
multiplication
\be
m=(m_{\cal C}\otimes m_{\cal D})(\id \otimes s \otimes \id),
\ee
$m_{\cal C}$ and $m_{\cal D}$ being the multiplication maps in
${\cal C}$ and ${\cal D}$,
where $s\colon {\cal D}\otimes {\cal C}\to {\cal C}\otimes {\cal
D}$ is a linear map satisfying
\be\label{crosp}
(\id\otimes s)(m_{\cal D}\otimes \id )=
(\id\otimes m_{\cal D})(s\otimes \id )(\id\otimes s),\qquad
(s\otimes\id )(\id\otimes m_{\cal C})=
(m_{\cal C}\otimes \id )(\id\otimes s)(s\otimes \id )
\ee
(this condition is equivalent to the associativity of $m$). For
unital algebras we require additionally $s(I\otimes c)=c\otimes
I$, $s(d\otimes I)=(I\otimes d)$
and for $^*$-algebras, we require that
$*_{12}*_{12}=\id$, where
\be
*_{12}=s(*\otimes *)\tau .
\ee
Under these conditions ${\cal C}\otimes {\cal D}$ becomes a
unital $^*$-algebra (called the {\em crossed tensor product of}
${\cal C}$ and ${\cal D}$) and the inclusions $c\mapsto c\otimes I$,
$d\mapsto I\otimes d$ are unital $^*$-homomorphisms (cf. for
instance \cite{DK}).
\begin{Thm}\label{Th}
If there exists a crossed tensor product of $*$-algebras $\Poly
(H)$ and $\Poly (V)$,  compatible with (\ref{cross}),
i.e. such that
\be
s(x^{A\cB}\otimes u^C_D)=T^{A\cB C}_{EK\cL}u^E_D\otimes x^{K\cL},
\ee
then it is unique. It exists if and only if
\be\label{S}
S=q^{-\frac12}M\qquad \mbox{or}\qquad S = q^{\frac12}M^{-1}
\ee
(square roots defined up to sign).
\end{Thm}
The proof is given in the Appendix (Section~\ref{calB}).

{}From now on we shall consider the case when $S=q^{-\frac12}M$
(the second case in (\ref{S}) is completely analogous).
We can write (\ref{xh}) as
\be
x_{\bf 1}h_{\bf 2}=\hat{W} h_{\bf 1}x_{\bf 2},
\ee
where
\be\label{WR}
\hat{W}=\hat{R}_- \;\;\;\mbox{for}\;\;\; \ov{q}=q,\qquad
\hat{W}=q^{-1}\hat{R}_-\;\;\; \mbox{for} \;\;\; |q|=1.
\ee

\section{Poincar\'{e} group with braided translations --- only
for $|q|=1$}

Now we can return to the problem when $\Dr$ preserves
(\ref{Mink}), i.e. when $P^{(-)}x_{\bf 1}x_{\bf 2} =0$ implies
$P^{(-)}\Dr x_{\bf 1}\Dr x_{\bf 2} =0$.
Assuming (\ref{Mink}), first two terms in
\be\label{DD}
\Dr x_{\bf 1}\Dr x_{\bf 2} = (x_{\bf 1}+h_{\bf 1}x'_{\bf
1})(x_{\bf 2}+h_{\bf 2}x'_{\bf 2})=x_{\bf 1}x_{\bf 2} +h_{\bf
1}x'_{\bf 1}h_{\bf 2}x'_{\bf 2} + x_{\bf 1}h_{\bf 2}x'_{\bf 2} +
h_{\bf 1}x'_{\bf 1}x_{\bf 2}
\ee
are obviously annihilated by $P^{(-)}$ (second, because
$P^{(-)}h_{\bf 1}h_{\bf 2}x'_{\bf 1}x'_{\bf 2}= h_{\bf 1}h_{\bf
2}P^{(-)}x'_{\bf 1}x'_{\bf 2}=0$). In the last term we shall
need to commute $x'_{\bf 1}$ with $x_{\bf 2}$. Normally they
just commute, but, it will be convenient to consider here the
following more general situation
\be
x_{\bf 1}'x_{\bf 2}=\hat{B}x_{\bf 1}x'_{\bf 2}\qquad \mbox{or}\qquad
x'_{\bf 2}x_{\bf 1}=Bx_{\bf 1}x'_{\bf 2}\;\;\; (\hat{B}=\tau B)
\ee
(for some matrix $B$). In particular, if $B=\id$, $x'_{\bf 1}$
and $x_{\bf 2}$ commute. Note that this more general assumption
does not affect previous results on the preservation of
(\ref{uuE}) and (\ref{cross}).

The sum of the last two terms in (\ref{DD}) is equal
$$
(\hat{W}h_{\bf 1}x_{\bf 2}x'_{\bf 1} + h_{\bf 1}x'_{\bf 1}x_{\bf
2})^{jk}=\hat{W}^{jk}_{ab}h^a_cx^bx'^c+ h^j_lB^{kl}_{bc}x^bx'^c=
(\hat{W}^{jk}_{ab}\dr ^l_c+\dr ^j_aB^{kl}_{bc})h^a_lx^bx'^c,
$$
hence, finally, $\Dr $ preserves (\ref{Mink}) when
\be\label{pwb}
P^{(-)}_{\bf 12}(\hat{W}_{\bf 12}+B_{\bf 23})=0.
\ee
Now, if $B=I$, then using (\ref{WR}) we see that the above
equality is possible only for $q^2=1$.

This is one more manifestation of the fact that the standard
$q$-deformation is not compatible with inhomogeneous groups.

On the other hand, if we could manage that
\be\label{main}
P^{(-)}(\hat{W}+\sigma I)=0\qquad\mbox{for some}\;\;\sigma ,
\ee
then $B=\sigma I$ satisfies (\ref{pwb}). In this case $\Dr $
preserves (\ref{Mink}) provided we consider the `braiding'
\be\label{braidx}
x'^jx^k=\sigma x^kx'^j.
\ee
Taking into account that $P^{(-)}$ is a projection and a
function of $\hat{W}$, condition (\ref{main}) means that
$P^{(-)}$ is a spectral projection of $\hat{W}$ corresponding to
a single eigenvalue (equal to $-\sigma$).
{}From (\ref{WR}), (\ref{R-1}) and (\ref{R-2}) it is clear that this is
possible only for $|q|=1$ and in this case $\sigma =q^{-1}$.

It is easy to check that (\ref{braidx}) define consistently a
crossed tensor product of ${\cal B}$ with itself.
Concluding, we have a family of braided Poincar\'{e} groups,
labelled by two parameters: $|q|=1$ and $t>0$.

\section{Minkowski space}

We present here explicitly the defining relations (\ref{Mink})
for the quantum Minkowski space corresponding to $|q|=1$ and $t>0$:
\begin{eqnarray*}
\ar\br & = & tq\br\ar \\
\ar\gr & = & t^{-1}q\gr\ar \\
\br\dr & = & tq\dr\br \\
\gr\dr & = & t^{-1}q\dr\gr \\
\br\gr & = & \gr\br \\
{} [\ar ,\dr ] & = & t^{-1}(q-q^{-1})\br\gr
\end{eqnarray*}
and
\be
\ar ^* =\ar,\qquad \dr ^*=\dr,\qquad \br ^*=\gr
\ee
(cf. (\ref{rel1})--(\ref{rel4}) and (\ref{5'})).
We have denoted the elements $x^{A\ov{B}}$ as follows
\be
x=\left( \ba{cc} x^{1\cj} & x^{1\cd }\\
        x^{2\cj } & x^{2\cd} \ea\right) =
  \left( \ba{cc} \ar & \br \\
                 \gr & \dr \ea\right) .
\ee
We may introduce a complex parameter $z:=q/t \neq 0$. The
corresponding quantum Minkowski space ${\cal M}_z$ is described
by the $*$-algebra $\Poly ({\cal M}_z)$ generated by elements
$\ar$, $\br$, $\gr$ satisfying
\be
\ar ^* =\ar,\qquad \dr ^*=\dr,\qquad \gr^*\gr =\gr\gr ^*,
\ee
such that
\begin{eqnarray*}
\ar\gr & = & z\gr\ar \\
\gr\dr & = & z\dr\gr \\
{} [\ar ,\dr ] & = & (z-\ov{z})\gr ^* \gr .
\end{eqnarray*}
The invariant {\em Minkowski length} (obtained as $E'_{12}(\tau
\ov{E'})_{34}X^{-1}_{23} x_{12}x_{34}$) is a central element of
$\Poly (V)$, equal
$$\frac{\ar\dr}{2z}+\frac{\dr\ar}{2\ov{z}}-\gr ^*\gr.$$

\section{Appendix}

\subsection{Braid intertwiners for Minkowski}\label{biM}

When we represent (\ref{Wpm}) on a diagram composed of elementary
crossings $X$, $X^{-1}$, $M$ and $K^{\pm 1}$, then it becomes
clear that it is sufficient to prove the `elementary moves'
\begin{eqnarray}
X^{-1}_{12}M_{23}X_{12} & = & X_{23}M_{12}X^{-1}_{23} \label{M1}\\
X^{-1}_{12}X^{-1}_{23}K^{\pm 1}_{12} & = &
K^{\pm 1}_{23}X^{-1}_{12}X^{-1}_{23} \label{M2} \\
K^{\pm 1}_{12}X_{23}X_{12} & = & X_{23}X_{12}K^{\pm 1}_{23}
\label{M3} \\
M_{12}X^{-1}_{23}X^{-1}_{12} & = & X^{-1}_{23}X^{-1}_{12}M_{23}
\label{M4} \\
X_{12}X_{23}M_{12} & = & M_{23}X_{12}X_{23} \label{M5} \\
X_{12}K^{\pm }_{23}X^{-1}_{12} & = & X^{-1}_{23}K^{\pm
1}_{12}X_{23} \label{M6}.
\end{eqnarray}
Equalities (\ref{M1}), (\ref{M4}) and (\ref{M5}) are mutually
equivalent. Also, (\ref{M2}), (\ref{M3}) and (\ref{M6}) are
mutually equivalent. Note that  (\ref{M3}) with `plus sign' is
obtained by taking the complex conjugation of (\ref{M5}). The
`minus sign' case is obtained by the complex conjugation of
$X_{12}X_{23}M^{-1}_{12}  =  M^{-1}_{23}X_{12}X_{23} $, which is
of course a simple consequence of (\ref{M5}) since $M^{-1}$ is a
polynomial of $M$.

Thus, it is sufficient to prove (\ref{M5}). This equality is
almost evident from \cite{WZ} (it expresses the fact that $X$
provides a representation of the standard $q$-commutation
relations). Let us prove it in detail.  Since $M$ is a linear
combination of $I$ and $EE'$, it is sufficient to show that
$$ X_{12}X_{23}E_{12}E'_{12} = E_{23}E'_{23}X_{12}X_{23} .$$
{}From (5) of \cite{WZ} we know that $X_{12}X_{23}E_{12}=cE_{23}$,
and, analogously, $E'_{23}X_{12}X_{23}=d E'_{12}$ for some
non-zero factors $c,d\in \bC$. But $d=c$, since
$$ (d E'_{12})E_{12} =(E'_{23}X_{12}X_{23})E_{12}=
E'_{23}(X_{12}X_{23})E_{12})=E'_{23}(cE_{23}),$$
and this ends the proof.

\subsection{Selection of parameters for Minkowski}\label{spM}

We shall write relations (\ref{Mink}) in explicit form.
The two cases of $X$ may be written in one formula
$$
X =e^1_{\ov{1}}\otimes e^{\ov{1}}_1+
e^2_{\ov{2}}\otimes e^{\ov{2}}_2+
t^{-1}(e^1_{\ov{2}}\otimes e^{\ov{2}}_1+e^2_{\ov{1}}\otimes
e^{\ov{1}}_2) +\er e^2_{\ov{1}}\otimes e^{\ov{2}}_1 ,
$$
where $\er =0$ in case 1 and $\er =\pm 1$, $t=q$ in case 2 (we
have rescaled $X$ for convenience). We have
$$
X^{-1} =e_1^{\ov{1}}\otimes e_{\ov{1}}^1+
e_2^{\ov{2}}\otimes e_{\ov{2}}^2+
t(e_1^{\ov{2}}\otimes e_{\ov{2}}^1+e_2^{\ov{1}}\otimes
e_{\ov{1}}^2) -\er e_1^{\ov{2}}\otimes e_{\ov{1}}^2 .
$$

Of course, (\ref{Mink}) is equivalent to
$$
(\cP '\otimes \cQ)X^{-1}_{23}x_{12}x_{34}=0\qquad
\mbox{and}\qquad
(\cP\otimes \cQ ' )X^{-1}_{23}x_{12}x_{34}=0.
$$
Using
$$\ker \cP '= \lel E\rr =\ker \lel E \rr ^\circ =\ker \lel
e^{11}, e^{22}, e^{21}+qe^{12}\rr ,$$
$$\ker \cQ  = \ker \tau \ov{E'} =\ker (\ov{q}e^{\ov{1}\ov{2}} -
e^{\ov{2}\ov{1}}),\qquad \ker \cP = \ker E' = \ker (qe^{21}-e^{12}),
$$
$$\ker \cQ '= \lel \tau \ov{E}\rr =\ker \lel \tau \ov{E} \rr
^\circ =\ker \lel e^{\ov{1}\ov{1}}, e^{\ov{2}\ov{2}},
\ov{q}e^{\ov{2}\ov{1}}+e^{\ov{1}\ov{2}}\rr ,$$
we see that
$\ker \cP '\otimes \cQ $ is composed of vectors which are
annihilated by the following three functionals
$$\{ e^{11},e^{22}, e^{21}+qe^{12}\}\otimes (\ov{q}e^{\ov{1}\ov{2}} -
e^{\ov{2}\ov{1}})=
$$
$$
= \{ \ov{q}e^{11\ov{1}\ov{2}}-e^{11\ov{2}\ov{1}},
\ov{q}e^{22\ov{1}\ov{2}}-e^{22\ov{2}\ov{1}},
|q|^2e^{12\ov{1}\ov{2}}+\ov{q}e^{21\ov{1}\ov{2}}-qe^{12\ov{2}\ov{1}}-
e^{21\ov{2}\ov{1}} \}
$$
(here $e^{11\ov{1}\ov{2}}:=e^{11}\otimes e^{\ov{1}\ov{2}}$ etc.)
and $\ker \cP\otimes \cQ ' $ is composed of vectors which are
annihilated by the following three functionals
$$  (qe^{21}-e^{12})\otimes \{e^{\ov{1}\ov{1}},e^{\ov{2}\ov{2}},
\ov{q}e^{\ov{2}\ov{1}}+e^{\ov{1}\ov{2}} \} =
$$
$$
=\{
qe^{21\ov{1}\ov{1}}-e^{12\ov{1}\ov{1}},qe^{21\ov{2}\ov{2}}-
e^{12\ov{2}\ov{2}},|q|^2e^{21\ov{2}\ov{1}}+qe^{21\ov{1}\ov{2}}-\ov{q}
e^{12\ov{2}\ov{1}}-e^{12\ov{1}\ov{2}} \} .
$$
Composing all the six functionals with $X^{-1}_{23}$, we get the
following functionals
$$
\{ \ov{q}(e^{1\ov{1}1\ov{2}} -\er e^{1\ov{2}2\ov{2}})-
te^{1\ov{2}1\ov{1}},
\ov{q}te^{2\ov{1}2\ov{2}}-e^{2\ov{2}2\ov{1}},
|q|^2te^{1\ov{1}2\ov{2}}+\ov{q}(e^{2\ov{1}1\ov{2}}- \er
e^{2\ov{2}2\ov{2}} )-qe^{1\ov{2}2\ov{1}} -te^{2\ov{2}1\ov{1}},
$$
$$
q(e^{2\ov{1}1\ov{1}} -\er
e^{2\ov{2}2\ov{1}})-te^{1\ov{1}2\ov{1}},
qte^{2\ov{2}1\ov{2}}-e^{1\ov{2}2\ov{2}},
|q|^2te^{2\ov{2}1\ov{1}}+q(e^{2\ov{1}1\ov{2}}- \er
e^{2\ov{2}2\ov{2}} )-\ov{q}e^{1\ov{2}2\ov{1}} -te^{1\ov{1}2\ov{2}}
\}
$$
and (\ref{Mink}) is equivalent to vanishing of these functionals
on $x_{12}x_{34}$. This way we obtain six following relations:
\begin{eqnarray*}
\ov{q}(x^{1\ov{1}}x^{1\ov{2}}-\er x^{1\ov{2}}x^{2\ov{2}}) -t
x^{1\ov{2}}x^{1\ov{1}} & = & 0 \\
\ov{q}tx^{2\ov{1}}x^{2\ov{2}} -x^{2\ov{2}}x^{2\ov{1}} & = & 0 \\
|q|^2tx^{1\ov{1}}x^{2\ov{2}} + \ov{q}(x^{2\ov{1}}x^{1\ov{2}}
-\er x^{2\ov{2}}x^{2\ov{2}}) -q
x^{1\ov{2}}x^{2\ov{1}}-tx^{2\ov{2}} x^{1\ov{1}} & = & 0 \\
q(x^{2\ov{1}}x^{1\ov{1}}-\er x^{2\ov{2}}x^{2\ov{1}}) -t
x^{1\ov{1}}x^{2\ov{1}} & = & 0 \\
qtx^{2\ov{2}}x^{1\ov{2}} -x^{1\ov{2}}x^{2\ov{2}} & = & 0 \\
|q|^2tx^{2\ov{2}}x^{1\ov{1}} + q(x^{2\ov{1}}x^{1\ov{2}}
-\er x^{2\ov{2}}x^{2\ov{2}}) -\ov{q}
x^{1\ov{2}}x^{2\ov{1}}-tx^{1\ov{1}} x^{2\ov{2}} & = & 0 .
\end{eqnarray*}
Substituting
\be
\left( \ba{cc} x^{1\cj} & x^{1\cd }\\
        x^{2\cj } & x^{2\cd} \ea\right) =
  \left( \ba{cc} \ar & \br \\
                 \gr & \dr \ea\right) ,
\ee
we can write these relations as follows
\begin{eqnarray}
\ov{q}(\ar \br -\er \br \dr ) -t
\br \ar  & = & 0 \label{ogo1}\\
\ov{q}t\gr \dr  -\dr \gr  & = & 0 \\
|q|^2t\ar \dr  + \ov{q}(\gr \br
-\er \dr \dr ) -q
\br \gr -t\dr  \ar  & = & 0 \\
q(\gr \ar -\er \dr \gr ) -t
\ar \gr  & = & 0 \\
qt\dr \br  -\br \dr  & = & 0 \\
|q|^2t\dr \ar  + q(\gr \br
-\er \dr \dr ) -\ov{q}
\br \gr -t\ar  \dr  & = & 0 \label{ogo6}.
\end{eqnarray}
Now we shall show that for the PBW theorem, condition
(\ref{warun}) is necessary. We thus consider the case 1 i.e.
$\er =0$. In this case, the commutation relations take the form
\begin{eqnarray}
\br\ar & = & \ov{q}t^{-1}\ar\br \label{rel1}\\
\gr\ar & = & q^{-1}t\ar\gr \label{rel2}\\
\dr\gr & = & \ov{q}t\gr\dr \label{rel3}\\
\dr\br & = & q^{-1}t^{-1}\br\dr \label{rel4}\\
|q|^2t\dr\ar + q\gr\br & = & t\ar\dr + \ov{q}\br\gr \label{rel5}\\
-t\dr\ar + \ov{q}\gr\br & = & -|q|^2t\ar\dr + q\br\gr \label{rel6}.
\end{eqnarray}
Taking $\ov{q}(\ref{rel5}) -q(\ref{rel6})$ and
$(\ref{rel5})+|q|^2(\ref{rel6})$ instead of (\ref{rel5}) and
(\ref{rel6}), we obtain:
\begin{eqnarray}
q(\ov{q}^2+1)t\dr\ar & = &
\ov{q}(q^2+1)t\ar\dr+(\ov{q}^2-q^2)\br\gr \label{re5} \\
q(\ov{q}^2+1)\gr\br & = & t(1-|q|^4)\ar\dr + \ov{q}(q^2+1)\br\gr
\label{re6}.
\end{eqnarray}
Using (\ref{rel1})--(\ref{rel4}) and (\ref{re5})--(\ref{re6}) it
is easy to see that each element of the algebra can be written
as a sum of (alphabetically) ordered monomials in
$\ar,\br,\gr,\dr$. Now, if we perform the two independent ways
of ordering of $q(\ov{q}^2+1)\gr\br\ar$, we obtain
$$
q(\ov{q}^2+1)\gr (\br\ar )=q(\ov{q}^2+1)\ov{q}t^{-1}\gr\ar\br =
\ov{q}(\ov{q}^2+1)\ar\gr\br = \frac{\ov{q}}{q}\ar
[t(1-|q|^4)\ar\dr +\ov{q}(q^2+1)\br\gr ]$$
on one hand, and
$$
q(\ov{q}^2+1)(\gr\br )\ar = [t(1-|q|^4)\ar\dr
+\ov{q}(q^2+1)\br\gr ]\ar =
\frac{1-|q|^4}{q(\ov{q}^2+1)}\ar [t\ov{q}(q^2+1)\ar\dr
+(\ov{q}^2-q^2)\br\gr ]+
$$
$$
+\frac{\ov{q}}{q}t(q^2+1)\br\ar\gr =
(1-|q|^4)\frac{\ov{q}}{q}\frac{q^2+1}{\ov{q}^2+1}t\ar\ar\dr +
\frac{(1-|q|^4)(\ov{q}^2-q^2)}{q(\ov{q}^2+1)}\ar\br\gr +
\frac{\ov{q}^2}{q}(q^2+1)\ar\br\gr
$$
on the other. Comparing the coefficients at $\ar\ar\dr$ we get
(\ref{warun}). Comparing at $\ar\br\gr$ gives exactly the same.
(We assume of course that $\ar\ar\dr$ and
$\ar\br\gr$ are linearly independent.)

For $|q|=1$, relations (\ref{re5}) and (\ref{re6}) are
equivalent to
\be\label{5'}
 \gr\br = \br\gr,\qquad [\ar ,\dr ] =
\frac{1}{t}(q-q^{-1})\br\gr ,
\ee
and the algebra $\Poly (V)$ resembles the usual $GL_q(2)$
algebra (it is the same, if $t=1$). One can easily check that
$\Poly (V)$ is a $q$-enveloping algebra in the sense of
\cite{B}, if we order the generators as follows:
$$e_1:=\ar ,\;\;\;\; e_2:=\br,\;\;\;\; e_3:=\gr,\;\;\;\;
e_4:=\dr ,$$
hence the PBW theorem holds in this case (see Theorem~2.8.1 of
\cite{B}).

If $\ov{q} =q$,  relations (\ref{re5}) and (\ref{re6}) are
equivalent to
$$ \dr\ar =\ar\dr,\qquad [\br ,\gr ] =t(q-q^{-1})\ar\dr .
$$
Replacing $\ar \leftrightarrow \br$, $\gr\leftrightarrow \dr$,
$t\leftrightarrow t^{-1}$, we get the same relations as in the
previous case, hence the PBW theorem holds also in this case.
Similarly, it holds for $\ov{q} = -q$.

The case when $\er =\pm 1$ and $q=t>0$ corresponds to the
standard quantum deformation of the Lorentz group, containing
as a subgroup $SU_q(2)$ or $SU_q(1,1)$ (depending on the sign of
$\er (q-1)$). Relations (\ref{ogo1})--(\ref{ogo6}) are then
equivalent to those considered by many authors
\cite{CW1,CW2,OZ}.
It can be easily shown that the PBW theorem holds in this case,
using the Diamond Lemma \cite{dl} (choose $\br <\ar <\dr <\gr $
as the total ordering).

\subsection{The algebra ${\cal B}$}\label{calB}

We set ${\cal A}:=\Poly (H)$, ${\cal C}:=\Poly (V)$

The uniqueness of $s$ is obvious, since its value on any
monomial can be reduced by (\ref{crosp}) to the case
(\ref{cross}).

Writing (\ref{cross}) as
$$ x_{12}u_3 = Tu_1x_{23}$$
(in case the crossed product exists), we get
\be
x_{12}E_{34}=x_{12}u_3u_4E_{34}=T_{123}u_1x_{23}u_4E_{14}=
T_{123}T_{234}u_1u_2x_{34}E_{12}=T_{123}T_{234}E_{12}x_{34},
\ee
hence $T$ must satisfy
\be\label{TTE}
T_{123}T_{234}E_{12} = E_{34}.
\ee
Taking into account that $X_{23}X_{12}E_{23}=E_{12}$, it means that
\be\label{SSE}
S_{12}S_{23}E_{12}=E_{23}.
\ee
It is easy to see that the only solutions of (\ref{SSE}) which
are intertwiners of $u_1u_2$ (hence of the form $aI+bEE'$) are
(\ref{S}).

Conversely, we shall show that if $S$ is given by (\ref{S}) then
there exists $s\colon \pC\otimes \pA\to \pA\otimes \pC$ with
required properties. Let $\tpA$ ($\tpC$) be the free
$^*$-algebra generated by $u^A_B$ ($x^{A\ov{B}}$). We have
\be
\pA = \tpA / \pJ _{\pA},\qquad \pC = \tpC / \pJ _{\pC},
\ee
where $\pJ _{\pA} = \lel \pJ ^0_{\pA}\rr$ is the ideal generated
by $\pJ ^0_{\pA} :=\{ u_1u_2E-E, E'u_1u_2-E',
Xu_1\ov{u}_2-\ov{u}_1u_2X \}$ in $\tpA$ and
$\pJ _{\pC} = \lel \pJ ^0_{\pC}\rr$ is the ideal generated
by $\pJ ^0_{\pC} :=\{ \hat{R}x_{\bf 1}x_{\bf 2} -x_{\bf 1}x_{\bf
2} , (x^{A\ov{B}})^*=x^{B\ov{A}} \}$ in $\tpC$. Here
$\hat{R}=\hat{R}_-$ for $\ov{q}=q$ and $\hat{R}=\hat{R}_+$ for
$|q|=1$ (cf. the discussion near (\ref{YBa})).
It is easy to see that there exists a (unique) map $\tis\colon
\tpC\otimes \tpA\to \tpA\otimes \tpC$ satisfying
(\ref{cross}) and (\ref{crosp}), with $\pA ,\pC, s$ replaced by
$\tpA,\tpC,\tis$.

The proof will be finished if we show that
$$
\tis (\tpC\otimes \pJ _{\pA})\subset \pJ _{\pA}\otimes \tpC,\qquad
\tis (\pJ _{\pC}\otimes \tpA)\subset \tpA\otimes \pJ _{\pC}.
$$
Since $\{ a\in \tpA : \tis (\tpC\otimes a)\subset
\pJ_{\pA}\otimes \tpC\}$ is an ideal in $\tpA$, it is sufficient
to show that
\be
\tis (\tpC\otimes \pJ ^0_\pA )\subset \pJ _{\pA}\otimes \tpC ,
\ee
and, similarly,
\be
\tis (\pJ ^0_{\pC}\otimes \tpA )\subset \tpA\otimes \pJ _{\pC}.
\ee
We shall show that
\be\label{codopo}
\tis (\tpC ^{(1)}\otimes \pJ ^0_\pA )\subset \pJ _{\pA}\otimes
\tpC ^{(1)} ,\qquad \tis (\pJ ^0_{\pC}\otimes \tpA ^{(1)} )\subset
\tpA ^{(1)}\otimes \pJ _{\pC},
\ee
where $\tpA ^{(1)}$ and $\tpC ^{(1)}$ denote the linear
subspaces spanned by the corresponding generators. It is
sufficient, because then from (\ref{crosp}) it follows that
\be
\tis (\tpC ^{(n)}\otimes \pJ ^0_\pA )\subset \pJ _{\pA}\otimes
\tpC ^{(n)} ,\qquad \tis (\pJ ^0_{\pC}\otimes \tpA ^{(n)} )\subset
\tpA ^{(n)}\otimes \pJ _{\pC},
\ee
where $\tpA ^{(n)}$ and $\tpC ^{(n)}$ denote the subspaces
spanned by monomials of order $n$.

To show (\ref{codopo}), note that
$$ x_{12}(u_3u_4E_{34}-E_{34}) =
T_{123}T_{234}u_1u_2x_{34}E_{12}-x_{12}E_{34} =
$$
$$
=T_{123}T_{234}(u_1u_2E_{12}-E_{12})x_{34}+T_{123}T_{234}E_{12}x_{34}
-x_{12}E_{34} = T_{123}T_{234}(u_1u_2E_{12}-E_{12})x_{34}
$$
belongs to $\pJ ^0_{\pA}\otimes \tpC$.
Similarly, $x_{12}(E'_{34}u_3u_4-E'_{34})\in \pJ ^0_{\pA}\otimes \tpC $
and
$$
x_{12}(X_{34}u_3\ov{u}_4 -\ov{u}_3u_4)=
X_{34}T_{123}T'_{234}u_1\ov{u}_2x_{34} -
T'_{123}T_{234}\ov{u}_1u_2x_{34}X_{12} = 0,$$
where $T'=  (\tau \ov{S}^{-1}\tau )_{23}X^{-1}_{12}$ is the
matrix appearing in (\ref{xcu}).
The equality $X_{34}T_{123}T'_{234}=
T'_{123}T_{234}X_{12}$ is proved using formulas of type
(\ref{M1})--(\ref{M6}).
Furthermore we have
$$
(\hat{R}_{1234}x_{12}x_{34}-x_{12}x_{34})u_5=
\hat{R}_{1234}x_{12}T_{345}u_3x_{45}-x_{12}T_{345}u_3x_{45}=
$$
$$
=\hat{R}_{1234}T_{345}T_{123}u_1x_{23}x_{45}-
T_{345}T_{123}u_1x_{23}x_{45}=
T_{345}T_{123}u_1(\hat{R}_{2345}x_{23}x_{45}-x_{23}x_{45})\in
\tpA\otimes \pJ ^0_{\pC},
$$
since
$\hat{R}_{1234}T_{345}T_{123}=T_{345}T_{123}\hat{R}_{2345}$ (it
also follows from (\ref{M1})--(\ref{M6})).

\section*{{\sl Acknowledgments}}

The author is very much indebted to Prof. S.L. Woronowicz for
drawing attention to the problem and valuable discussions.

\end{document}